\documentclass[10pt]{article}
\usepackage{float}
\usepackage{graphicx}
\usepackage{rotating}
\usepackage{caption}
\usepackage{subcaption}
\usepackage[linesnumbered]{algorithm2e}
\usepackage{comment}
\usepackage{multicol}
\usepackage{url}
\usepackage[margin=0.5in]{geometry}

\title{A Progressive Statistical Method for Preconditioning Matrix-Free Solution of High-Order Discretization of Linear Time-Dependent Problems}
\author{\textit{A. Ghasemi} and \textit{L. K. Taylor} \\ SimCenter: National Center for Comp. Eng., \\ 701 East M L King Blvd, Chattanooga, TN 37403}

\date{}

\begin{document}

\begin{multicols*}{2}
\maketitle

\section*{Introduction and Motivations}
Consider the general form of linear time-dependent problems
\begin{equation}
\label{eq_gen_lin_time_dep}
\frac{\partial u}{\partial t} = R = \mathcal{L}\left(t,u\right) = \sum_{k}\sum_{i} \alpha_{ki} \frac{\partial^k u}{{\partial x_i }^k}, 
\end{equation}
where the residual vector $R$ contains high-order derivatives of the dependent variable $u(t,x_1,x_2,\ldots,x_n)$ with respect to independent variables. Equation (\ref{eq_gen_lin_time_dep}) is an abstract representation of a generic Partial Differential Equation (PDE) in the residual form where the linear operator $\mathcal{L}$ accounts for a linear combination of the k$^{th}$ partial derivative in the i$^{th}$ direction, i.e. $\partial^k /{\partial x_i }^k$. When discretized using a suitable discretization method like Finite-Differences, Finite Elements or Finite Volume \cite{FEM}, the resulting equations can be written in the following semi-discrete form 
\begin{equation}
\label{eq_gen_lin_time_dep_semi_disc_1}
\frac{d }{dt}\mathbf{u} = \mathbf{L}\;\mathbf{u}.
\end{equation}
The fully discrete from is obtained after eq.(\ref{eq_gen_lin_time_dep_semi_disc_1}) is discretized in time. There are many choices available for time-discretization \cite{Butcher}. However a common step in these algorithms, which is the most costly part of the solution, is simply an Euler implicit scheme presented below
\begin{equation}
\label{eq_gen_lin_time_dep_semi_disc_2}
\frac{\mathbf{u}^{n+1} - \mathbf{u}^{n}}{\Delta t} = \mathbf{L} \mathbf{u}^{n+1},
\end{equation}
which leads to the solution of the following \textit{system of linear equations}
\begin{equation}
\label{eq_gen_lin_time_dep_semi_disc_3}
\underbrace{\left(\mathbf{I} - \Delta t \mathbf{L}\right)}_{\mathbf{A}}\mathbf{u}^{n+1} = \mathbf{u}^{n},
\end{equation}  

where $\mathbf{A} = \left(\mathbf{I} - \Delta t \mathbf{L}\right)$ is typically a huge sparse matrix for practical problems. For this reason, a direct solution of eq.(\ref{eq_gen_lin_time_dep_semi_disc_3}) is impossible and thus an iterative algorithm is usually used for the solution. These algorithms range from stationary methods including Jacobi, Gauss-Seidel, SOR to more sophisticated Krylov subspace methods CG (Conjugate Gradient) and GMRES (Generalized Minimal Residual). The reader should refer to \cite{templates} and references therein for more detail. In this work, the GMRES algorithm \cite{Saad_1986}\footnote{The least-squares is used in Ref.\cite{Saad_1986} to solve a system containing an upper Hessenberg matrix. However in this work it is used to estimate the preconditioner from a regression point of view.} is used as the base algorithm however the statistical preconditioning proposed here can be readily applied to all Krylov subspace methods because they are based on matrix-vector multiplication.  

The important point is that for many practical problems the matrix $\mathbf{A}$ in eq.(\ref{eq_gen_lin_time_dep_semi_disc_3}) is so huge that it does not fit in the physical memory of the current machines. This situation typically happens in Large Eddy Simulation (LES) or Direct Numerical Simulations (DNS) of a turbulent fluid flow \cite{Lesieur}. In this situation a full implicit solution procedure is avoided to prevent the problem of storing the full Jacobian matrix. Yet still another possibility is available by using a matrix-free implementation \cite{Knoll}. In this approach the matrix $\mathbf{A}$ is never computed/stored but eq.(\ref{eq_gen_lin_time_dep_semi_disc_3}) is solved by replacing the matrix-vector product in the GMRES algorithm with residual computation. However the only disadvantage of this approach is that preconditioning is almost impossible since entries $a_{ij}$ of $\mathbf{A}$ are not available. This is the \textit{motivation} for current work where a simple statistical method is used to estimate $\mathbf{A}$ by least-squaring the history of matrix-vector product obtained in the GMRES algorithm.  \section*{The PSP-GMRES Algorithm} 
\begin{figure}[H]
  \centering
  \includegraphics[trim = 32mm 26mm 15mm 7mm, clip, width=0.5\textwidth]{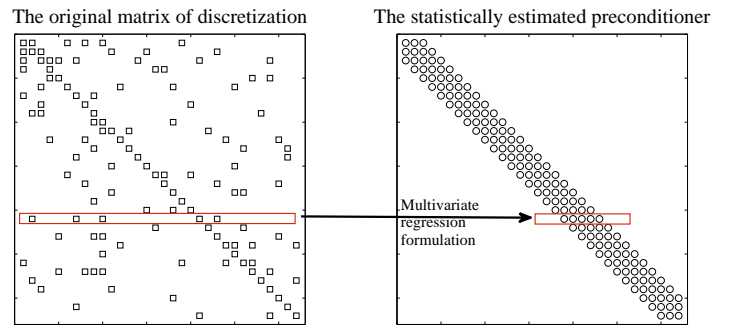}
  \caption{The schematic of the statistical reduction of the sparse matrix $\mathbf{A}$ into the banded diagonal matrix $\mathbf{N}$.}\label{fig_precond_sketch}
\end{figure}

The details of the Progressive Statistical Preconditioning for GMRES is covered in algorithms (\ref{alg_PSPGMRES}) and (\ref{alg_1MREP}) below. To solve the linear system eq.(\ref{eq_gen_lin_time_dep_semi_disc_3}), the original GMRES algorithm starts as usual in algorithm (\ref{alg_PSPGMRES}). However when a matrix-vector product is performed in lines 4 and 11 of (\ref{alg_PSPGMRES}), for the given vector $\mathbf{P}_x(:,i)$, the result is stored in the i$^{th}$ column of matrix $\mathbf{P}_y$, i.e. $\mathbf{P}_y(:,i)$. Therefore the \textit{statistical} dataset $\mathbf{P}_y$ versus $\mathbf{P}_x$ is obtained \textit{progressively} in a sense that the more matrix-vector product is performed the better dataset is obtained. Now, the primary goal is to find a relation between $\mathbf{P}_y$ and $\mathbf{P}_x$ such that it has minimum error compared to the exact relation $\mathbf{P}_y = \mathbf{A} \mathbf{P}_x$. This relation, which is estimated using a banded diagonal matrix $\mathbf{N}$ such that $\mathbf{P}_y \approx \mathbf{N} \mathbf{P}_x$ will be further used as a preconditioner in the GMRES algorithm (\ref{alg_PSPGMRES}) when eq.(\ref{eq_gen_lin_time_dep_semi_disc_3}) is solved for the next time step $\Delta t$. The approximation of $\mathbf{A}$ with the banded diagonal matrix $\mathbf{N}$ is illustrated in fig.(\ref{fig_precond_sketch}). As shown, all off-diagonal entries on a given rows of the original huge matrix $\mathbf{A}$ are reduced to a banded matrix $\mathbf{N}$ where the matrix-vector product would be cheap. This reduction procedure is completely described in algorithm (\ref{alg_1MREP}) using a multi-regression estimation of matrix $\mathbf{P}_y$ versus matrix $\mathbf{P}_x$. As shown, the output of the GMRES algorithm, i.e. $\mathbf{P}_x$ and $\mathbf{P}_y$ are sent to MREP (Multi regressor Preconditioner) where for $d=1$, a tridiagonal estimation (three regressors) is enforced. The emphasis on three regressor model is based on the fact that this approach leads to a tridiagonal matrix $\mathbf{N}$ which can be solved efficiently in $O(n)$ FLOPs using Thomas algorithm. Thus it is expected that this three-regressor model would save a lot when the preconditioned equations need to be solved in lines 10 and 42 of alg.(\ref{alg_PSPGMRES}).  

{\small 
\IncMargin{1em}
\begin{algorithm}[H]
  \SetAlgoLined
  {\bf subroutine}  (X, R, $\mathbf{P}_x$, $\mathbf{P}_y$) $\gets$ {\bf PSPGMRES} ($\mathbf{A}$, B, $\mathbf{N}$, $\epsilon$, $n_A$, $X_0$, $n_r$)\;
  $i \gets 1$\;
  \For{$l \leftarrow 1$ \KwTo $n_r$}
  { $\mathbf{P}_x(:,i) \gets X_0$, $\mathbf{P}_y(:,i) \gets \mathbf{A} \; X_0$\;
    $\bar{R} \gets B- \mathbf{P}_y(:,i)$, $i \gets i + 1$\;
    $\mathbf{V}(:,1) \gets \bar{R}/\|\bar{R}\|$ \;
    $G(1) \gets \|\bar{R}\|$\;
    \For{$k \leftarrow 1$ \KwTo $n_A$}
    {
      $G(k+1) \gets $0\;
      $Y_k \gets {\mathbf{N}}^{-1} \mathbf{V}(:,k)$\;
      $\mathbf{P}_x(:,i) \gets Y_k$, $\mathbf{P}_y(:,i) \gets \mathbf{A} \; Y_k$\;
      $U_k \gets \mathbf{P}_y(:,i)$, $i \gets i + 1$\;
      \For{$j \gets 1$ \KwTo $k$}
      {
        $\mathbf{H}(j,k) \gets {\mathbf{V}(:,j)}^{'} \; U_k$ \;
        $U_k \gets U_k - \mathbf{H}(j,k) \; \mathbf{V}(:,j)$ \;
      }
      
      $\mathbf{H}(k+1,k) \gets \|U_k\|$\;        
      $\mathbf{V}(:,k+1) \gets U_k/\mathbf{H}(k+1,k)$\;        
      \For{$j \gets 1$ \KwTo $k-1$}
      {
        $\delta \gets \mathbf{H}(j,k)$\;
        $\mathbf{H}(j,k) \gets \delta C(j) + S(j) \mathbf{H}(j+1,k)$\;
        $\mathbf{H}(j+1,k) \gets -\delta S(j) + C(j) \mathbf{H}(j+1,k)$\;
      }
      $\gamma \gets \sqrt{{\mathbf{H}(k,k)}^{2} + {\mathbf{H}(k+1,k)}^{2}}$\;
      $C(k) \gets \mathbf{H}(k,k) / \gamma$\;
      $S(k) \gets \mathbf{H}(k+1,k) / \gamma$\;
      $\mathbf{H}(k,k) \gets \gamma$\;
      $\mathbf{H}(k+1,k) \gets 0$\;
      $\delta \gets G(k)$\;
      $G(k) \gets C(k) \delta + S(k) G(k+1)$\;
      $G(k+1) \gets -S(k)  \delta + C(k) G(k+1)$\;
      $R(k) \gets \|G(k+1)\|$\;
      \If{ $R(k) \leq \epsilon$ }
      {
        \textit{finish-flag} $gets$ 1\;
        {\bf break}\;
      }

    }  
    $Q \gets {\mathbf{H}}^{-1} G$\;
    \For {$j \gets 1$ \KwTo $k$}
    {
      $Z_k \gets  Z_k + Q(j) \mathbf{V}(:,j)$\;
    }
    $X \gets X_0 + {\mathbf{N}}^{-1} Z_k$\;    
    $X_0 \gets X$\;    
    {\bf clean} $G, \mathbf{V}, \mathbf{H}, C, S$\;   
    \If{finish-flag is set} {
      {\bf break}\;
    }

  }
  \caption{The Progressive Statistical Preconditioned GMRES algorithm (PSPGMRES) with possible restarting.}\label{alg_PSPGMRES}
\end{algorithm}
\DecMargin{1em}
} As a summary, the output of the matrix-vector product in the GMRES alg.(\ref{alg_PSPGMRES}) is stored in two matrices $\mathbf{P}_x$ and $\mathbf{P}_y$ which are used as dataset for further modeling of the unavailable matrix $\mathbf{A}$. Subsequently, in the next immediate stage, the matrices $\mathbf{P}_x$ and $\mathbf{P}_y$ are used in alg. (\ref{alg_1MREP}) where a three regressors least square procedure (four unknowns) is used (d=1) to obtain the tridiagonal entries of the approximate matrix $\mathbf{N}$. This matrix is guaranteed to approximate the exact matrix $\mathbf{A}$ with minimal residual. Finally the computed matrix $\mathbf{N}$ is used as a preconditioner in the GMRES algorithm (\ref{alg_PSPGMRES}) to improve the convergence. Obviously, the more convergence of GMRES is improved the closer $\mathbf{N}$ to $\mathbf{A}$ is and vice versa. Now, it is interesting to see the results of this method in the implementation. {\small
\IncMargin{1em}
\begin{algorithm}[H]
  \SetAlgoLined
  {\bf function}  $\mathbf{N}$ $\gets$ {\bf MREP} ($\mathbf{P}_x$, $\mathbf{P}_y$, $d$)\;
  n $\gets$ number of rows of $\mathbf{P}_x$\;
  m $\gets$ number of columns of $\mathbf{P}_x$\;

  \For{$i \gets 1$ \KwTo $d$}
  {
    ($\beta_0$, $\beta_1$) $\gets$ {\bf linear fit}($\mathbf{P}_x(i,:)$, $\mathbf{P}_y(i,:)$)\;
    $\mathbf{N}(i,i) \gets \beta_1$\;
  }

  \For{$i \gets (d+1)$ \KwTo $(n-d)$}
  {
    $\mathbf{P}_x^{*} \gets \left[ \begin{array}{c}
        {\left[\begin{array}{cccc} 1 & 1 & \ldots & 1 \end{array}\right]}_{1\times m} \\
        \mathbf{P}_x(i-d,:) \\
        \vdots \\ 
        \mathbf{P}_x(i-1,:) \\
        \mathbf{P}_x(i,:) \\ 
        \mathbf{P}_x(i+1,:)\\
        \vdots\\
        \mathbf{P}_x(i+d,:)
      \end{array} \right] $\;
    $\mathbf{P}_y^{*} \gets \mathbf{P}_y(i,:)$\;
    $N^{*} \gets {\left(\mathbf{P}_x^{*}(\mathbf{P}_x^{*})^{'}\right)}^{-1} \mathbf{P}_x^{*} \left({\mathbf{P}_y^{*}}^{'}\right)$\;
    
    \For{$j \gets 2$ \KwTo $2(d+1)$}
    {
      $\mathbf{N}(i,i-j-d+1) \gets N^{*}(j)$\; 
    }
  }

  \For{$i \gets (n-d+1)$ \KwTo $n$}
  {
    ($\beta_0$, $\beta_1$) $\gets$ {\bf linear fit}($\mathbf{P}_x(i,:)$, $\mathbf{P}_y(i,:)$)\;
    $\mathbf{N}(i,i) \gets \beta_1$\;
  }

  \caption{The algorithm for Multi-Regressor Estimator of the Preconditioner (MREP).}\label{alg_1MREP}
\end{algorithm}
\DecMargin{1em}
}
\section*{Results}
\begin{figure}[H]
  \centering
  \includegraphics[trim = 2mm 2mm 6mm 1mm, clip, width=0.24\textwidth]{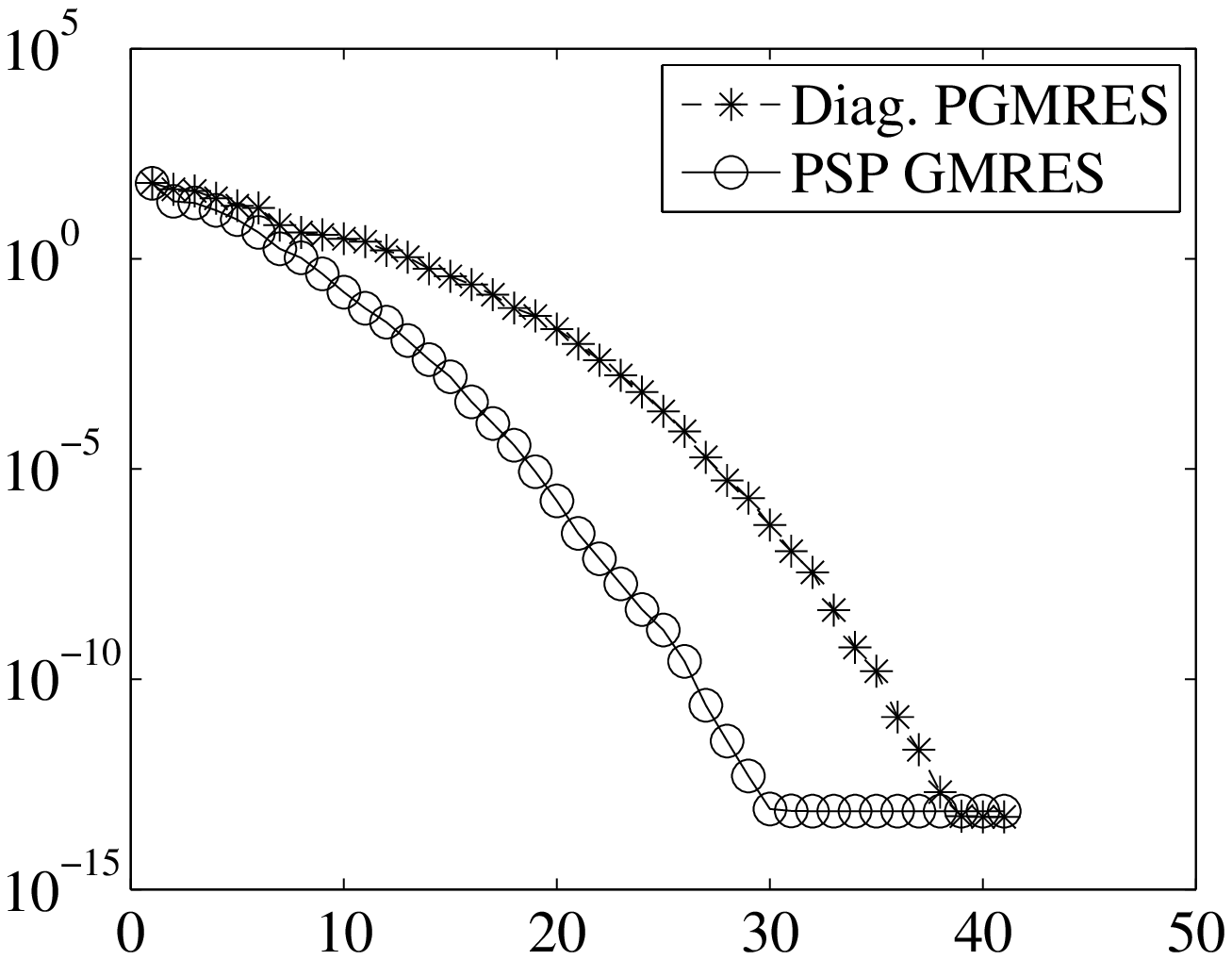}
  \includegraphics[trim = 2mm 2mm 6mm 1mm, clip, width=0.24\textwidth]{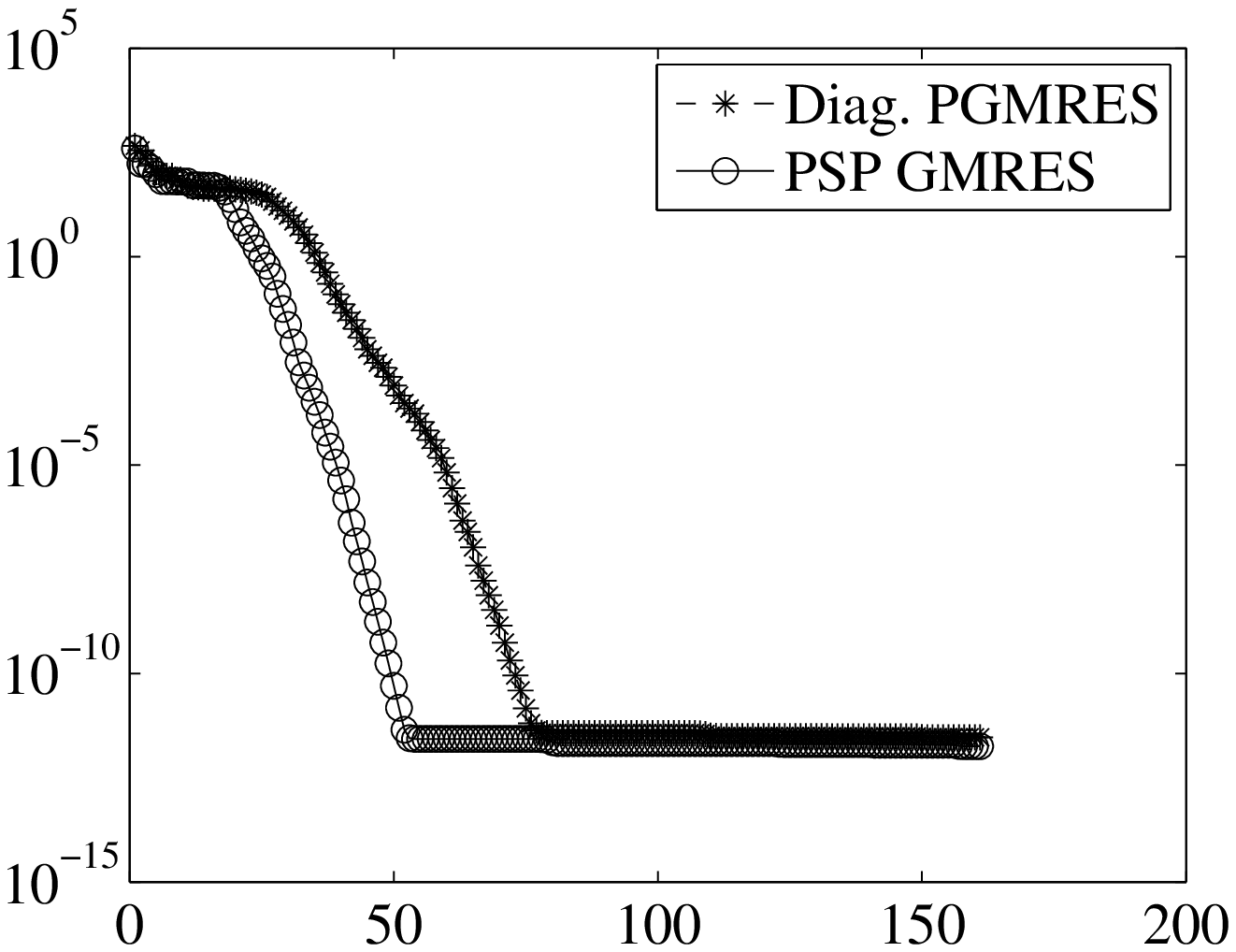}
  \includegraphics[trim = 2mm 2mm 6mm 1mm, clip, width=0.24\textwidth]{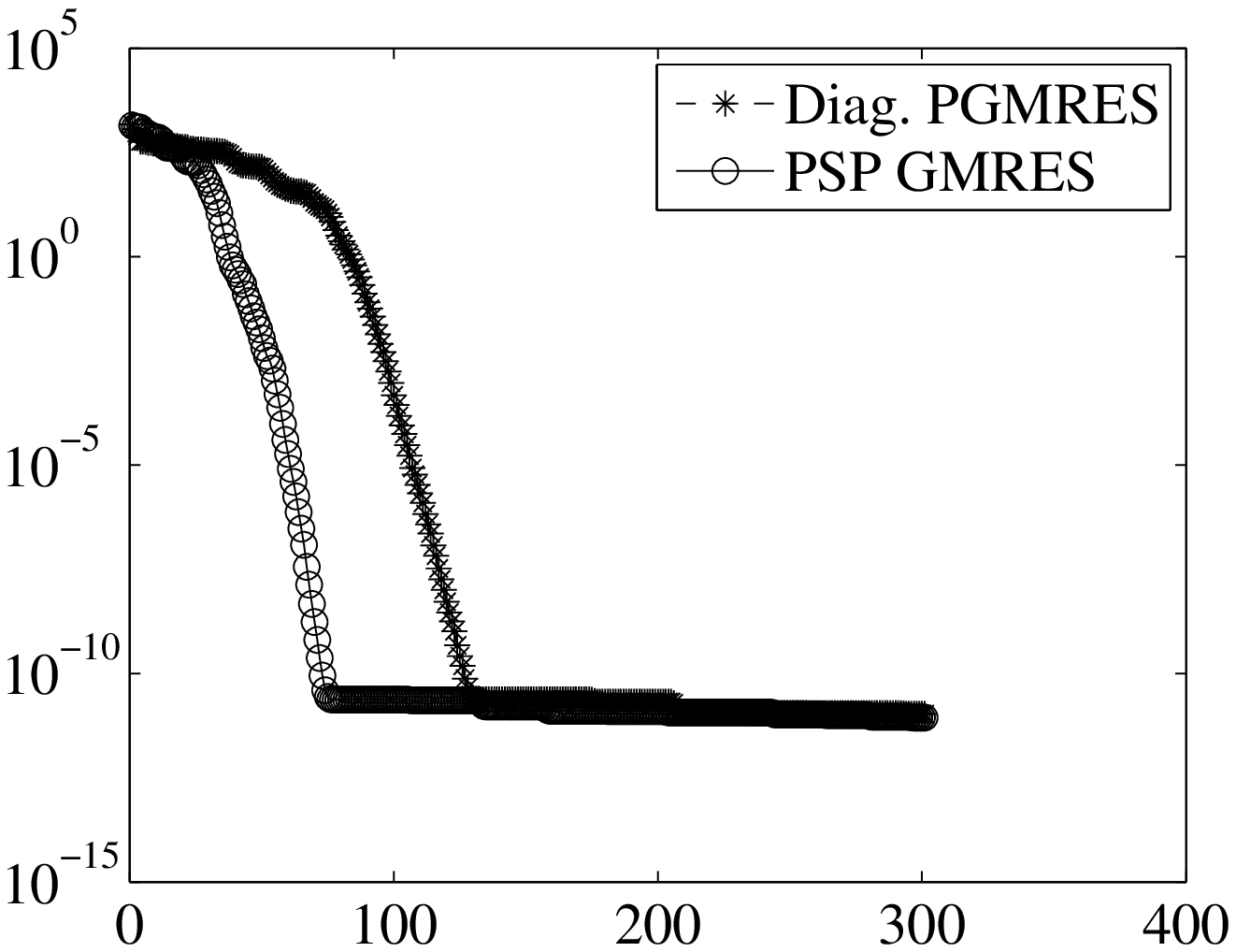}
  \includegraphics[trim = 2mm 2mm 6mm 1mm, clip, width=0.24\textwidth]{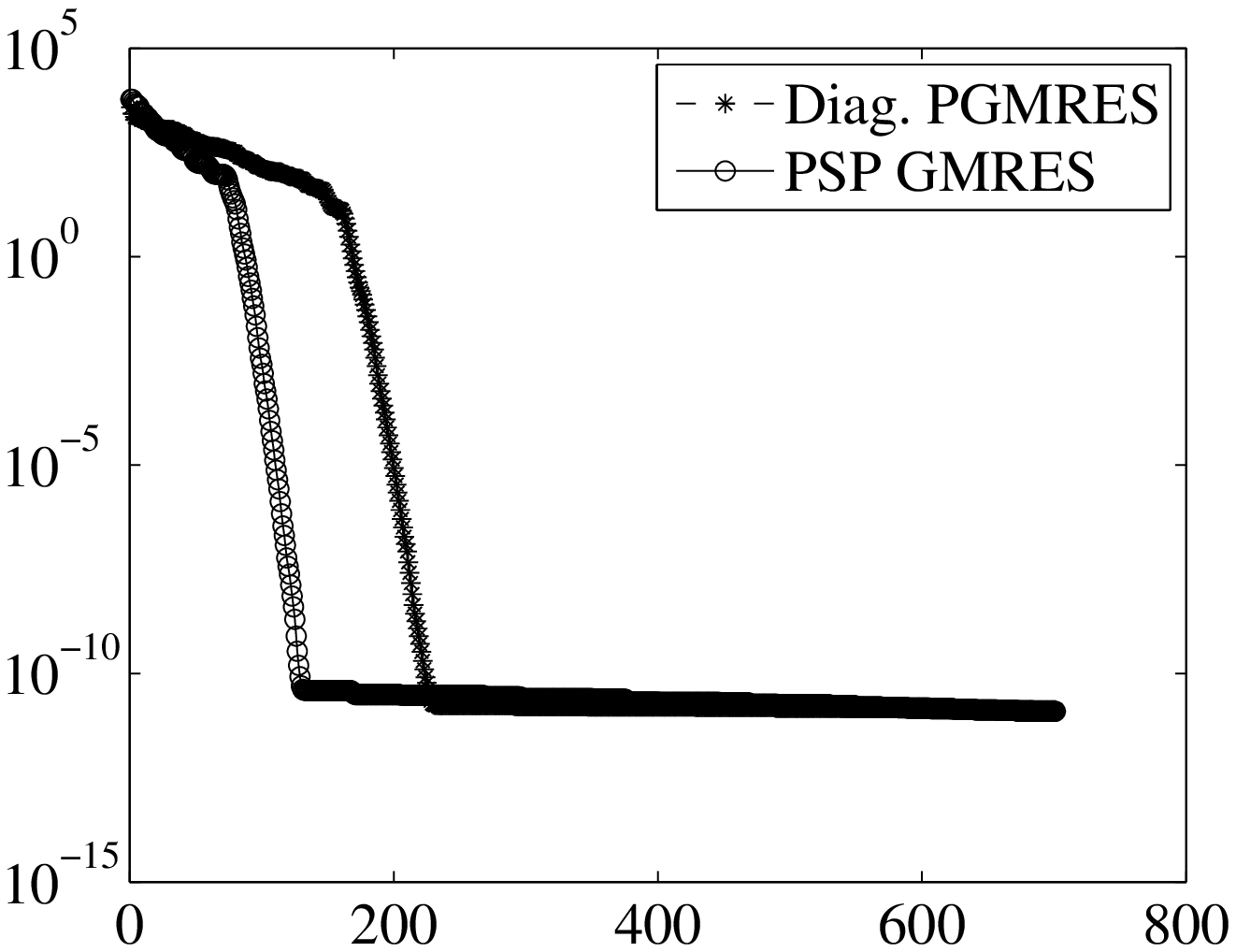}
  \includegraphics[trim = 2mm 2mm 6mm 1mm, clip, width=0.4\textwidth]{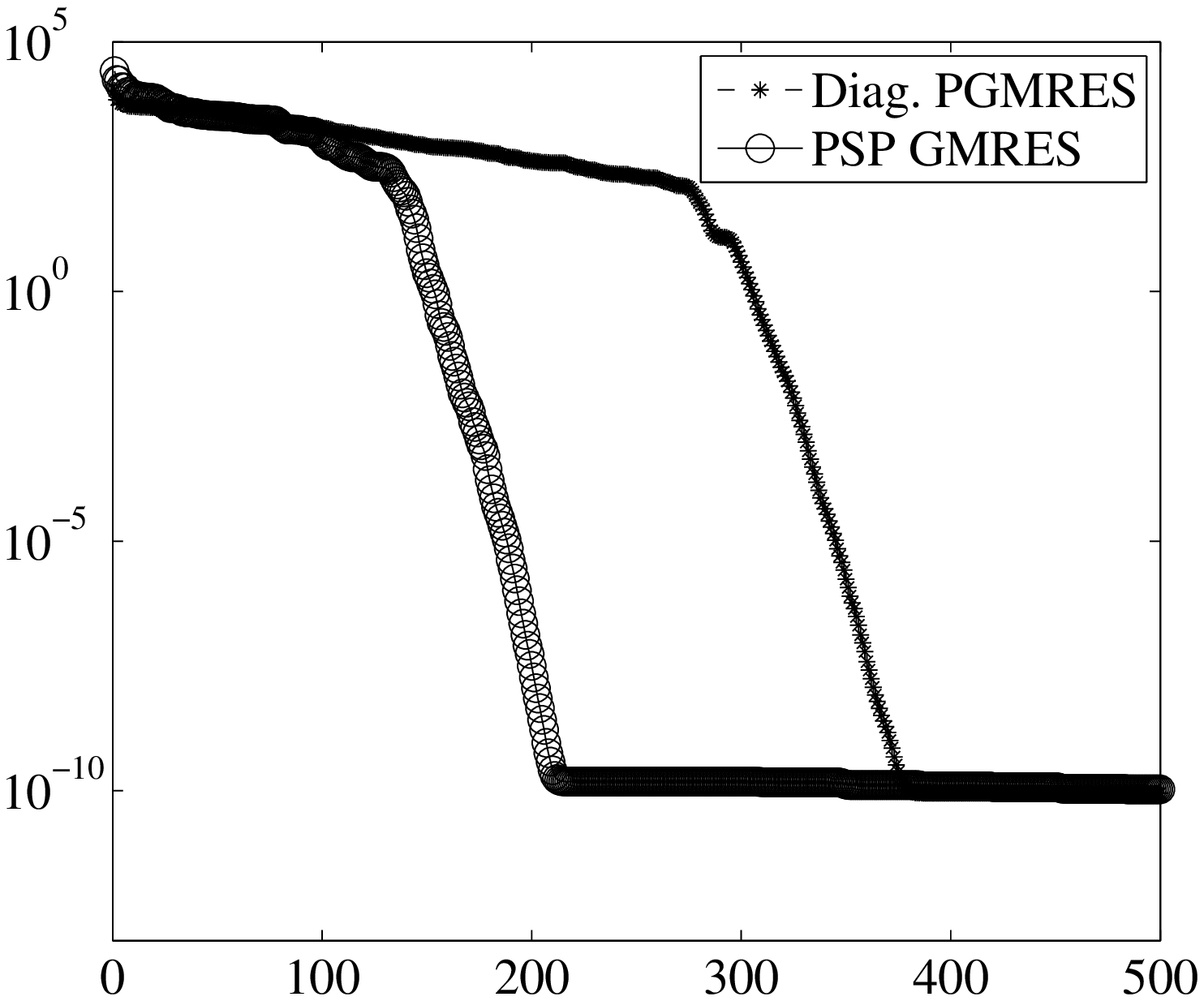}
  \caption{The residual norm versus number of GMRES iterations. From top-left to bottom-right: 20x20, 80x80, 150x150, 350x350 and 700x700.}\label{fig_conv_comp}
\end{figure}

\begin{figure}[H]
  \centering
  \includegraphics[trim = 28mm 7mm 22mm 7mm, clip, width=0.24\textwidth]{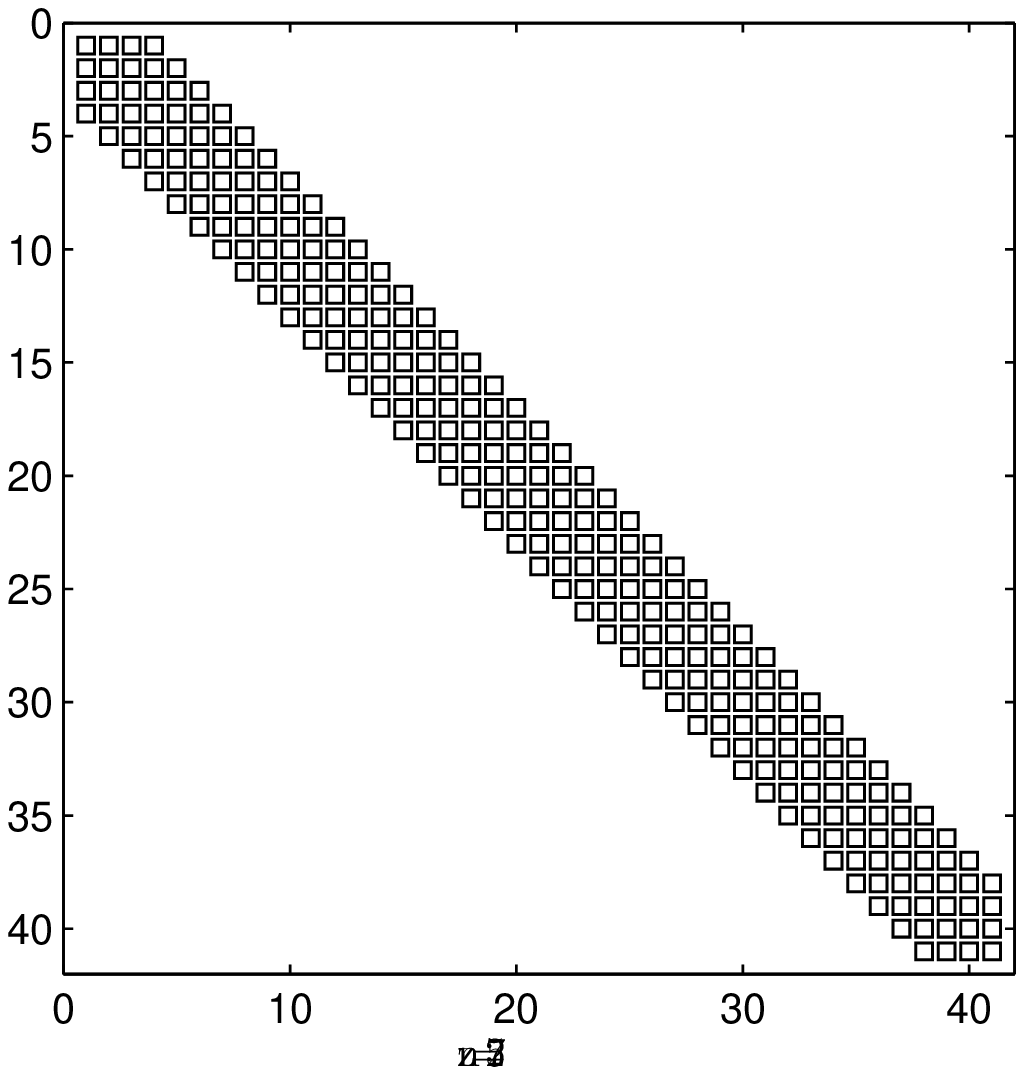}
  \includegraphics[trim = 28mm 7mm 22mm 7mm, clip, width=0.24\textwidth]{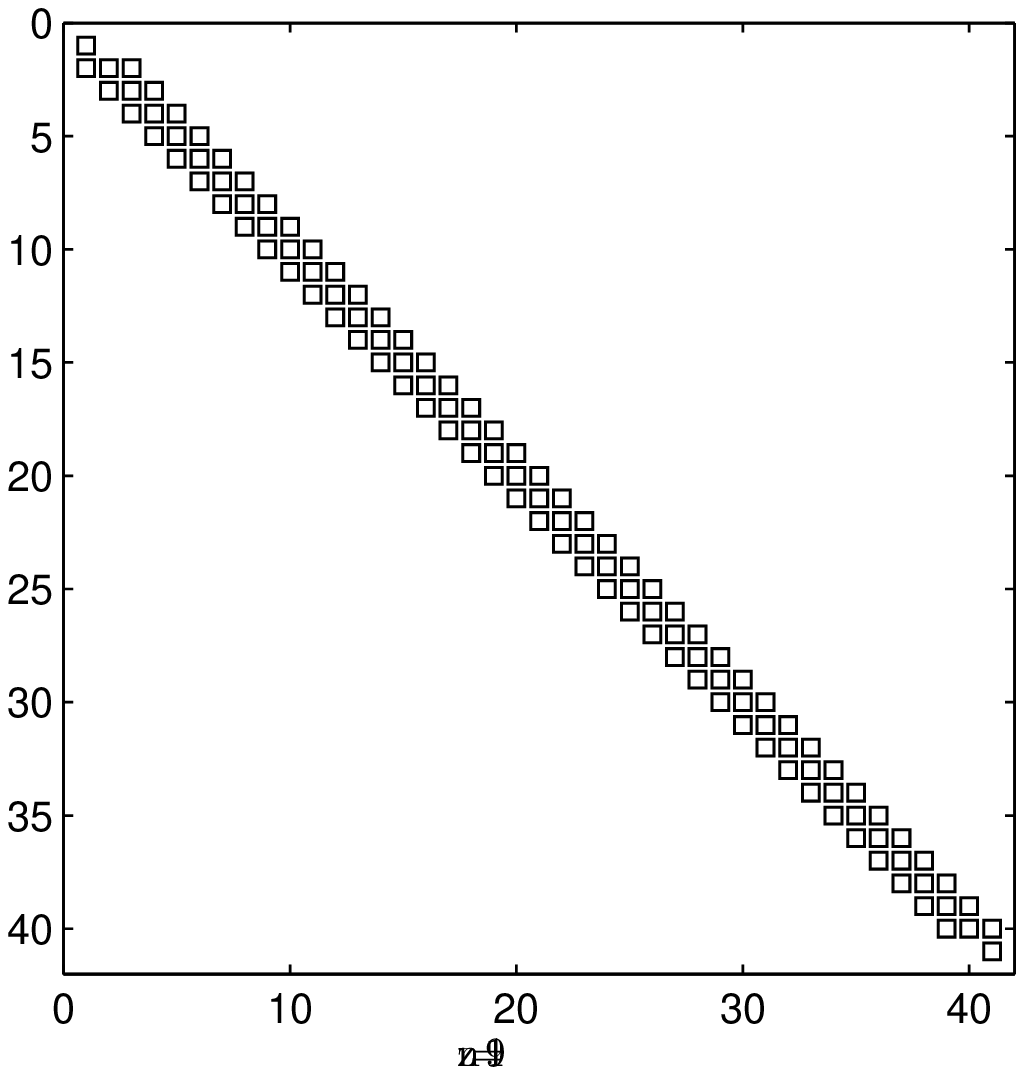}
  \caption{Left) The original pattern of the matrix. Right) The statistically reduced matrix used as preconditioner. The size of the matrix is $n=20$.}\label{fig_mat_comp}
\end{figure}

The algorithms are implemented in a computer program. In a mock test case, the matrix $\mathbf{A}_{n\times n}$ in eq.(\ref{eq_gen_lin_time_dep_semi_disc_3}) is a seven diagonal random matrix which is generated on fly in a way that diagonal dominance is preserved. This vector $\mathbf{u}^{n}$ is just a vector from 1 to $n$. The size of the system of equation, i.e. $n$ is varied to $n=\{20,80,150,350, 700\}$. The structure of the original matrix $\mathbf{A}$ for the case 20x20 and the preconditioner $\mathbf{N}$ is shown in fig.(\ref{fig_mat_comp}). As shown for $d=1$, a tridiagonal matrix is generated. Figure (\ref{fig_conv_comp}) compares the convergence of the PSP-GMRES method with the original GMRES algorithm.
As shown for case 20x20, the convergence is decreased from 40 iterations to 30 iterations. This indicates a factor 4/3 in the speed-up of original algorithm. Interestingly more impressive results would be obtained for larger systems. According to the same figure, as the size of the linear system is increased, the speed-up genrated by the PSP algorithm is also increased. As shown in fig.(\ref{fig_conv_comp}), the new statistical preconditioning approach reduces the convergence of the original GMRES from approximately 400 iterations to approximately 200 iterations when $\mathbf{A}$ is 700x700 matrix! This is a numerically validated indication of the robustness of the statistical approach.  
 \section*{More Regressors and Conclusions} It should be noted that the number of regressors used in alg.(\ref{alg_1MREP}) was chosen to be three ($d=1$) in the presented implementation in order to obtain a tridiagonal preconditioner. However other alternatives might be available for $d>1$ for future works. For example, all off-diagonal entries might be selected as regressors in the initial model. Therefore the initial model might have $n$ regressors. As the GMRES progress and more matrix-vector multiplications are performed, a standard model selection procedure, like backward elimination, might be used to eliminate inappropriate regressors\cite{Book,Gao}. This would result in optimum model which might result in optimum speed-up of the convergence of the GMRES algorithm. However this approach has two \textit{expensive} parts that should be analyzed mathematically in a separate paper. 

The first disadvanage is that elimination of inappropriate regressors requires to generate some models at startup. This would require \textit{at least} the cost of $O(n^3)$ for generating the first model using all regressors which is not acceptable because the cost is even higher than solving the system \textit{without} preconditioning! In addition, more models are needed to be created and compared to each other during this procedure which makes it even more expensive.

The second disadvangate of increasing the number of regressors is that the final preconditioner matrix would not be tridiagonal therefore the cost of solution of preconditioned system in lines 10 and 42 of alg.(\ref{alg_PSPGMRES}) would increase dramatically in this case! 
As a final conclusion, the three regressor model presented here has already generated impressive results in improving the convergence of the original GMRES algorithm according to fig(\ref{fig_conv_comp}). However, there might be situations that adding more regressors, i.e. $d>1$ would result in outstanding improvements in the convergence rate. It was mentioned that a backward elimination method is not practical in this case thus a forward selection is appropriate. This situation should be considered more rigorously in the future works.

{\tiny

}

\end{multicols*}

\end{document}